\documentclass{amsproc}
\usepackage{breqn,float,amssymb,pdflscape,rotating,graphicx}
\makeatletter
\@namedef{subjclassname@2020}{%
  \textup{2020} Mathematics Subject Classification}
\makeatother
\newtheorem{theorem}{Theorem}[section]

\theoremstyle{definition}

\newtheorem{example}[theorem]{Example}

\theoremstyle{remark}

\numberwithin{equation}{section}
\allowdisplaybreaks
\begin{document}

\title{Extended Watson-Harkins Sum}

%    Remove any unused author tags.

%    author one information
\author{Robert Reynolds}
\address[Robert Reynolds]{Department of Mathematics and Statistics, York University, Toronto, ON, Canada, M3J1P3}
\email[Corresponding author]{milver@my.yorku.ca}
\thanks{}

%    author two information
%\author{ Allan Stauffer}
%\address[Allan Stauffer]{Department of Mathematics and Statistics, York University, Toronto, ON, Canada, M3J1P3}
%\email{stauffer@yorku.ca}
%\thanks{This research is supported by NSERC Canada under Grant 504070}



\date{}

\dedicatory{}

\begin{abstract}
The Watson-Harkins sum involving the product of the cosine and cosecant functions is extended to derive the finite sum of generalized Hurwitz-Lerch Zeta functions is derived in terms of the Hurwitz-Lerch Zeta function. A transformation formula arises for various finite values of the parameters involved. The finite product of trigonometric functions are also derived. All the results in this work are new.
\end{abstract}

\maketitle
\section{Introduction}
The work done on theories of the structure of the atom by Watson \cite{watson} involved the finite sum of the cosecant function which was used in the calculation of the angular velocity of the system. Harkins et al. \cite{harkins} also produced work on the structure of the atom where the finite sum of the product the cosine and cosecant functions were used to calculate the Nicholson period of vibration. Since the finite sum of trigonometric functions are of high usage and importance the authors embark on using their method to derive a finite sum of the generalized cases of trigonometric functions namely the Hurwitz-Lerch Zeta function. \\\\
In this present work we derive the finite sum given by
\begin{multline}
\sum_{p=0}^{n-1}e^{-\frac{2 i j (\pi  l p+m n)}{n}} \left(\Phi \left(e^{\frac{2 i (m n+l p \pi )}{n}},-k,-j-\frac{1}{2} i \log
   (a)\right)\right.\\
   \left.+e^{2 i (2 j+1) \left(\frac{\pi  l p}{n}+m\right)} \Phi \left(e^{\frac{2 i (m n+l p \pi
   )}{n}},-k,j-\frac{1}{2} i \log (a)+1\right)\right)\\\\
=   (2 i)^{-k} n \left(\log ^k(a)+2^{k+1} (i n)^k e^{2 i m n} \Phi
   \left(e^{2 i m n},-k,1-\frac{i \log (a)}{2 n}\right)\right)
\end{multline}
where the variables $k,a,m$ are general complex numbers and $j=0,1, \ldots , n-1$; $l$ and $n$ are relatively prime numbers. This new expression is then used to derive special cases in terms of trigonometric functions. \\\\
Some special cases evaluated are in terms of the log-gamma function and Glaisher's constant $A$ which is featured in Barnes’ $G$ -Function in section (5.17) in \cite{dlmf} and Bernoulli number in section (24.2) in \cite{dlmf}. New functional identities are also derived which are used in the study of the distribution of prime numbers see section (25.16(i)) in \cite{dlmf}, the evaluation of Euler sums see section (25.16(ii)) in \cite{dlmf}. According to section (25.17) in \cite{dlmf}, the zeta function occurs  in the evaluation of the partition function of ideal quantum gases (both Bose–Einstein and Fermi–Dirac cases), and it determines the critical gas temperature and density for the Bose–Einstein condensation phase transition in a dilute gas (Lifshitz and Pitaevski\u{i} \cite{Lifshitz}). Quantum field theory often encounters formally divergent sums that need to be evaluated by a process of regularization: for example, the energy of the electromagnetic vacuum in a confined space (Casimir–Polder effect). It has been found possible to perform such regularizations by equating the divergent sums to zeta functions and associated functions (Elizalde \cite{Elizalde}).\\\\
Our preliminaries start with the application of a contour integral method to a finite sum. Let $a$, $k$, $m$ and $w$ be general complex numbers and $n\in\mathbb{Z^{+}}$ and $l$ is a relative prime number, the contour integral form \cite{reyn4}, of the finite sum given by equation (4.4.6.12) in \cite{prud1} is given by
\begin{multline}\label{fs-c1}
\frac{1}{2\pi i}\int_{C}\sum_{p=0}^{n-1}a^w w^{-k-1} \csc \left(\frac{\pi  l p}{n}+m+w\right) \cos \left((2 j+1) \left(\frac{\pi  l
   p}{n}+m+w\right)\right)dw\\
   =\frac{1}{2\pi i}\int_{C}n a^w w^{-k-1} \cot (n (m+w))dw
\end{multline}
The derivations follow the method used by us in \cite{reyn4}. This method involves using a form of the generalized Cauchy's integral formula given by
\begin{equation}\label{intro:cauchy}
\frac{y^k}{\Gamma(k+1)}=\frac{1}{2\pi i}\int_{C}\frac{e^{wy}}{w^{k+1}}dw,
\end{equation}
where $y,w\in\mathbb{C}$ and $C$ is in general an open contour in the complex plane where the bilinear concomitant \cite{reyn4} has the same value at the end points of the contour. This method involves using a form of equation (\ref{intro:cauchy}) then multiplies both sides by a function, then takes the finite sum of both sides. This yields a finite sum in terms of a contour integral. Then we multiply both sides of equation (\ref{intro:cauchy})  by another function and take the infinite sum of both sides such that the contour integral of both equations are the same.
\section{The Hurwitz-Lerch Zeta Function}

We use equation (1.11.3) in \cite{erd} where $\Phi(z,s,v)$ is the Lerch function which is a generalization of the Hurwitz Zeta $\zeta(s,v)$ and Polylogarithm functions $Li_{n}(z)$. The Lerch function has a series representation given by

\begin{equation}\label{knuth:lerch}
\Phi(z,s,v)=\sum_{n=0}^{\infty}(v+n)^{-s}z^{n}
\end{equation}
where $|z|<1, v \neq 0,-1,-2,-3,..,$ and is continued analytically by its integral representation given by

\begin{equation}\label{knuth:lerch1}
\Phi(z,s,v)=\frac{1}{\Gamma(s)}\int_{0}^{\infty}\frac{t^{s-1}e^{-vt}}{1-ze^{-t}}dt=\frac{1}{\Gamma(s)}\int_{0}^{\infty}\frac{t^{s-1}e^{-(v-1)t}}{e^{t}-z}dt
\end{equation}
where $Re(v)>0$, and either $|z| \leq 1, z \neq 1, Re(s)>0$, or $z=1, Re(s)>1$.
\section{Contour Integral Representation For The Finite Sum Of The Hurwitz-Lerch Zeta Functions}
We use the method in \cite{reyn4}. The cut and contour are in the first quadrant of the complex $w$-plane with $0 < Re(w+m)$.  The cut approaches the origin from the interior of the first quadrant and goes to infinity vertically and the contour goes round the origin with zero radius and is on opposite sides of the cut. Using a generalization of Cauchy's integral formula (\ref{intro:cauchy}) we first replace $y$ by $ \log (a)+i x+y$ then multiply both sides by $e^{itx}$ then form a second equation by replacing $x$ by $-x$ and adding both equations to get
\begin{multline}
\frac{e^{-i t x} \left(e^{2 i t x} (\log (a)+i x+y)^k+(\log (a)-i x+y)^k\right)}{\Gamma(k+1)}\\
=\frac{1}{2\pi i}\int_{C}2 a^w w^{-k-1} e^{w y}
   \cos (x (t+w))dw
\end{multline}
Next we replace $y$ by $i(2y+1)$, $t$ by $\frac{\pi  l p}{n}+m$, $x$ by $2j+1$ and multiply both sides by $e^{i (2 y+1) \left(\frac{\pi  l p}{n}+m\right)}$ and take the infinite and finite sums over $y\in[0,\infty)$ and $p\in [0,n-1]$, respectively and simplify in terms of the Hurwitz-Lerch Zeta function to get
\begin{multline}\label{fsci}
-\frac{1}{\Gamma(k+1)}i^{k+1} 2^k \sum_{p=0}^{n-1}e^{-\frac{2 i j (\pi  l p+m n)}{n}} \left(\Phi \left(e^{\frac{2 i (m n+l p \pi
   )}{n}},-k,-j-\frac{1}{2} i \log (a)\right)\right.\\
   \left.+e^{2 i (2 j+1) \left(\frac{\pi  l p}{n}+m\right)} \Phi \left(e^{\frac{2
   i (m n+l p \pi )}{n}},-k,j-\frac{1}{2} i \log (a)+1\right)\right)\\\\
   =-\frac{1}{2\pi i}\sum_{y=0}^{\infty}\sum_{p=0}^{n-1}\int_{C}2 i a^w
   w^{-k-1} e^{\frac{i (2 y+1) (\pi  l p+n (m+w))}{n}} \cos \left((2 j+1) \left(\frac{\pi  l
   p}{n}+m+w\right)\right)dw\\\\
    =-\frac{1}{2\pi i}\sum_{p=0}^{n-1}\int_{C}\sum_{y=0}^{\infty}2 i a^w
   w^{-k-1} e^{\frac{i (2 y+1) (\pi  l p+n (m+w))}{n}} \cos \left((2 j+1) \left(\frac{\pi  l
   p}{n}+m+w\right)\right)dw\\\\
=\frac{1}{2\pi i}\int_{C} \sum_{p=0}^{n-1}  a^w w^{-k-1} \csc \left(\frac{\pi  l
   p}{n}+m+w\right) \cos \left((2 j+1) \left(\frac{\pi  l p}{n}+m+w\right)\right)dw\\\\
   =\frac{1}{2\pi i}\int_{C}n a^w w^{-k-1} \cot (n (m+w))dw
\end{multline}
from equation (4.4.6.12) in \cite{prud1} where $Im(n(m+w))>0$ and $Re(n(m+w))>0, j=0, \ldots ,n-1$; $l$ and $n$ are relatively prime numbers. We apply Tonelli's theorem for multiple sums, see page 177 in \cite{gelca} as the summands are of bounded measure over the space $\mathbb{C} \times [0,n] $.
\section{Contour Integral Representation For The Hurwitz-Lerch Zeta Function}
We use the method in \cite{reyn4}. Using equation (\ref{intro:cauchy})  we first replace $\log (a)+2 i n (y+1)$ and multiply both sides by $-2 i n e^{2 i m n (y+1)}$ then take the infinite sum over $y\in [0,\infty)$ and simplify in terms of the Hurwitz-Lerch Zeta function to get
\begin{multline}\label{isci}
-\frac{2^{k+1} (i n)^{k+1} e^{2 i m n} \Phi \left(e^{2 i m n},-k,1-\frac{i \log (a)}{2 n}\right)}{\Gamma(k+1)}\\\\
=-\frac{1}{2\pi i}\sum_{y=0}^{\infty}\int_{C}2 i n a^w
   w^{-k-1} e^{2 i n (y+1) (m+w)}dw\\\\
=-\frac{1}{2\pi i}\int_{C}\sum_{y=0}^{\infty}2 i n a^w
   w^{-k-1} e^{2 i n (y+1) (m+w)}dw\\\\
=\frac{1}{2\pi i}\int_{C}n a^w
   w^{-k-1} \cot (n (m+w))+i n a^w w^{-k-1}dw
\end{multline}
from equation (1.232.1) in \cite{grad} where the $Im(n(m+w))>0$ in order for the sum to converge.
\subsection{The additional contour}
We use the method in \cite{reyn4}. Using equation (\ref{intro:cauchy}) we replace $y$ by $\log(a)$ and multiply both sides by $ni$ and simplify to get
\begin{equation}\label{addc}
\frac{i n \log ^k(a)}{\Gamma(k+1)}=\frac{1}{2\pi i}\int_{C}i n a^w w^{-k-1}dw
\end{equation}
\section{The finite sum of Hurwitz-Lerch Zeta functions in terms of the Hurwitz-Lerch Zeta function}
\begin{theorem}
For all $k,a,m\in\mathbb{C}$, $l$ is prime, $j<n-1$ then,
\begin{multline}\label{fslf}
\sum_{p=0}^{n-1}e^{-\frac{2 i j (\pi  l p+m n)}{n}} \left(\Phi \left(e^{\frac{2 i (m n+l p \pi )}{n}},-k,-j-\frac{1}{2} i \log
   (a)\right)\right.\\
   \left.+e^{2 i (2 j+1) \left(\frac{\pi  l p}{n}+m\right)} \Phi \left(e^{\frac{2 i (m n+l p \pi
   )}{n}},-k,j-\frac{1}{2} i \log (a)+1\right)\right)\\\\
=   (2 i)^{-k} n \left(\log ^k(a)+2^{k+1} (i n)^k e^{2 i m n} \Phi
   \left(e^{2 i m n},-k,1-\frac{i \log (a)}{2 n}\right)\right)
\end{multline}
\end{theorem}
\begin{proof}
The right-hand sides of relation (\ref{fsci}), and the addition of relations (\ref{addc}) and (\ref{isci}) are identical; hence, the left-hand sides of the same are identical too. Simplifying with the Gamma function yields the desired conclusion.
\end{proof}
\begin{example}
The Degenerate Case.
\begin{equation}
\sum_{p=0}^{n-1}\csc \left(\frac{\pi  l p}{n}+m\right) \left(-\cos \left((2 j+1) \left(\frac{\pi  l
   p}{n}+m\right)\right)\right)=-n \cot (m n)
\end{equation}
\end{example}
\begin{proof}
Use Equation (\ref{fslf}) and set $k=0$ and simplify using entry (2) in Table below (64:12:7) in \cite{atlas}.
\end{proof}
\section{Trigonometric Product-Recursion Identities And Transformation Formula}
Finite trigonometric products are important in many areas of mathematics and their uses and properties are detailed in \cite{peter} page 201. In this section we will evaluate equation (\ref{fslf}) for various parameter values and derive new finite trigonometric products.
\begin{example}
A Finite Product Involving The Exponential Of Trigonometric Functions
\begin{multline}
\prod_{p=0}^{n-1}\sin ^3\left(\frac{\pi  p}{n}+\frac{x}{2}\right) \csc ^2\left(\frac{\pi  p}{n}+\frac{x}{4}\right) \csc
   \left(\frac{\pi  p}{n}+x\right)\\
    \exp \left(4 \sin ^2\left(\frac{x}{4}\right) \left(2 \cos \left(\frac{2 \pi 
   p}{n}+x\right)+\cos \left(\frac{2 \pi  p}{n}+\frac{3 x}{2}\right)\right)\right)\\\\
   =2 \cos ^2\left(\frac{n x}{4}\right)
   \sec \left(\frac{n x}{2}\right)
\end{multline}
\end{example}
\begin{proof}
Use equation (\ref{fslf}) and set $k=1,a=1,m=x,l=1,j=1$ and simplify using the method in section (8.1) in \cite{reyn_ejpam}.
\end{proof}
\begin{example}
A Finite Product Involving The Exponential Of Trigonometric Functions
\begin{multline}
\prod_{p=0}^{n-1}\left(\sin \left(\frac{\pi  p}{n}+\frac{x}{2}\right) \csc \left(\frac{\pi  p}{n}+x\right)\right)^{i \pi }\\
 \exp
   \left(\cot \left(\frac{\pi  p}{n}+\frac{x}{2}\right)-\cot \left(\frac{\pi  p}{n}+x\right)+2 \sin
   \left(\frac{x}{2}\right)\right.\\
   \left.
    \left(2 \cos \left(\frac{2 \pi  p}{n}+\frac{3 x}{2}\right)+i \pi  \sin \left(\frac{2 \pi 
   p}{n}+\frac{3 x}{2}\right)\right)\right)\\\\
=   2^{-i \pi } \sec ^{i \pi }\left(\frac{n x}{2}\right) (\sinh (n \csc (n
   x))+\cosh (n \csc (n x)))
\end{multline}
\end{example}
\begin{proof}
Use equation (\ref{fslf}) and set $k=1,a=-1,m=x,l=1,j=1$ and simplify using the method in section (8.1) in \cite{reyn_ejpam}.
\end{proof}
\begin{example}
A Finite Product Involving The Exponential Of Trigonometric Functions
\begin{multline}
\prod_{p=0}^{n-1}\sin ^3\left(\frac{\pi  p}{n}+\frac{x}{2}\right) \csc ^2\left(\frac{\pi  p}{n}+\frac{x}{4}\right) \csc
   \left(\frac{\pi  p}{n}+x\right)\\
    \exp \left(-4 \sin \left(\frac{x}{4}\right) \sin \left(\frac{2 \pi  p}{n}+\frac{3
   x}{4}\right)-\frac{4}{3} \sin \left(\frac{3 x}{4}\right)\right.\\
   \left. \sin \left(\frac{6 \pi  p}{n}+\frac{9 x}{4}\right)+\sin
   (x) \sin \left(\frac{4 \pi  p}{n}+3 x\right)+\frac{2}{3} \sin \left(\frac{3 x}{2}\right)\right.\\
   \left. \sin \left(\frac{6 \pi 
   p}{n}+\frac{9 x}{2}\right)-4 \sin \left(\frac{x}{2}\right) \sin \left(\frac{\pi  p}{n}\right) \cos \left(\frac{3
   \pi  p}{n}+\frac{3 x}{2}\right)\right)\\\\
=   \sec \left(\frac{n x}{2}\right)+1
\end{multline}
\end{example}
\begin{proof}
Use equation (\ref{fslf}) and set $k=1,a=1,m=x,l=1,j=3$ and simplify using the method in section (8.1) in \cite{reyn_ejpam}.
\end{proof}
\begin{example}
A Finite Product Involving The Exponential Of Trigonometric Functions
\begin{multline}
\prod_{p=0}^{n-1}\sin ^{14}\left(\frac{2 \pi  p}{n}+\frac{x}{4}\right) \sin \left(\frac{2 \pi  p}{n}+x\right) \csc
   ^8\left(\frac{2 \pi  p}{n}+\frac{x}{8}\right) \csc ^7\left(\frac{2 \pi  p}{n}+\frac{x}{2}\right)\\
    \exp \left(-8 \cos
   \left(\frac{4 \pi  p}{n}+\frac{x}{4}\right)+14 \cos \left(\frac{4 \pi  p}{n}+\frac{x}{2}\right)-7 \cos
   \left(\frac{4 \pi  p}{n}+x\right)\right.\\
   \left.+\cos \left(\frac{4 \pi  p}{n}+2 x\right)\right)\\\\
=   8 \cos ^8\left(\frac{n
   x}{8}\right) \cos \left(\frac{n x}{2}\right) \sec ^6\left(\frac{n x}{4}\right)
\end{multline}
\end{example}
\begin{proof}
Use equation (\ref{fslf}) and set $k=2,a=1,m=x,l=2,j=1$ and simplify using the method in section (8.1) in \cite{reyn_ejpam}.
\end{proof}
\begin{example}
A  Recurrence Identity With Consecutive Neighbours.
\begin{equation}
\Phi (z,s,a)=z^j \left(z^{j+1} (-\Phi (z,s,a+2 j+1))+2 z \Phi (z,s,a+j+1)+i^s (i (a+j))^{-s}\right)
\end{equation}
\end{example}
\begin{proof}
Use equation (\ref{fslf}) and set $n=1,m=\log(z)/i,k=-s,a=e^{ai},j=0$ then set $a\to 2 (a-1),z\to \sqrt{z}$ and simplify.
\end{proof}
\begin{example}
A  Recurrence Identity With Consecutive Neighbours.
\begin{equation}
\Phi (z,s,a)=z \Phi (z,s,a+1)+a^{-s}
\end{equation}
\end{example}
\begin{proof}
Use equation (\ref{fslf}) and set $n=1,m=\log(z)/i,k=-s,a=e^{ai},j=0$ and simplify. This is equation (10.06.17.0002.01) in \cite{wolfram}.
\end{proof}
\begin{example}
A  Recurrence Identity With Consecutive Neighbours in terms of Glaisher's constant $A$.
\begin{equation}
-\Phi'(-i,-1,0)-\Phi'(i,-1,0)-\text{Li}_{-1}'(- i)-\text{Li}_{-1}'(i)=\log \left(\frac{A^{24}}{16\ 2^{2/3} e^2}\right)
\end{equation}
\end{example}
\begin{proof}
Use equation (\ref{fslf}) and set $n=2,a=1,j=0,m=\pi/4,l=3$ and simplify. Then take the first partial derivative with respect to $k$ and set $k=1$ and simplify.
\end{proof}
\begin{example}
A Finite Product Of Quotient Tangent Functions In Terms Of The Exponential Of The Hurwitz-Lerch Zeta Function. Plots of the right-hand are produced for a special case.
\begin{multline}
\prod_{p=0}^{n-1}\left(i \cot \left(\frac{\pi  l p+m n}{2 n}\right)\right)^{2 e^{-\frac{i (\pi  l p+m n)}{n}}} \left(-i \tan \left(\frac{\pi  l p+n r}{2 n}\right)\right)^{2 e^{-\frac{i
   (\pi  l p+n r)}{n}}}\\
=\exp \left(2 e^{2 i m n} \Phi \left(e^{2 i m n},1,1+\frac{1}{2 n}\right)-2 e^{2 i n r} \Phi \left(e^{2 i n r},1,1+\frac{1}{2 n}\right)\right)
\end{multline}
\end{example}
\begin{proof}
Use equation (\ref{fslf}) and form a second equation by replacing $m \to r$ and take their difference. Next set $k=-1,a=e^{i},j=-1$ and simplify in terms of the arctangent function using entry (3) in Table below (64:12:7) in \cite{atlas}. Next take the exponential of both sides and simplify. 
\end{proof}
\begin{figure}[H]
\includegraphics[scale=0.5]{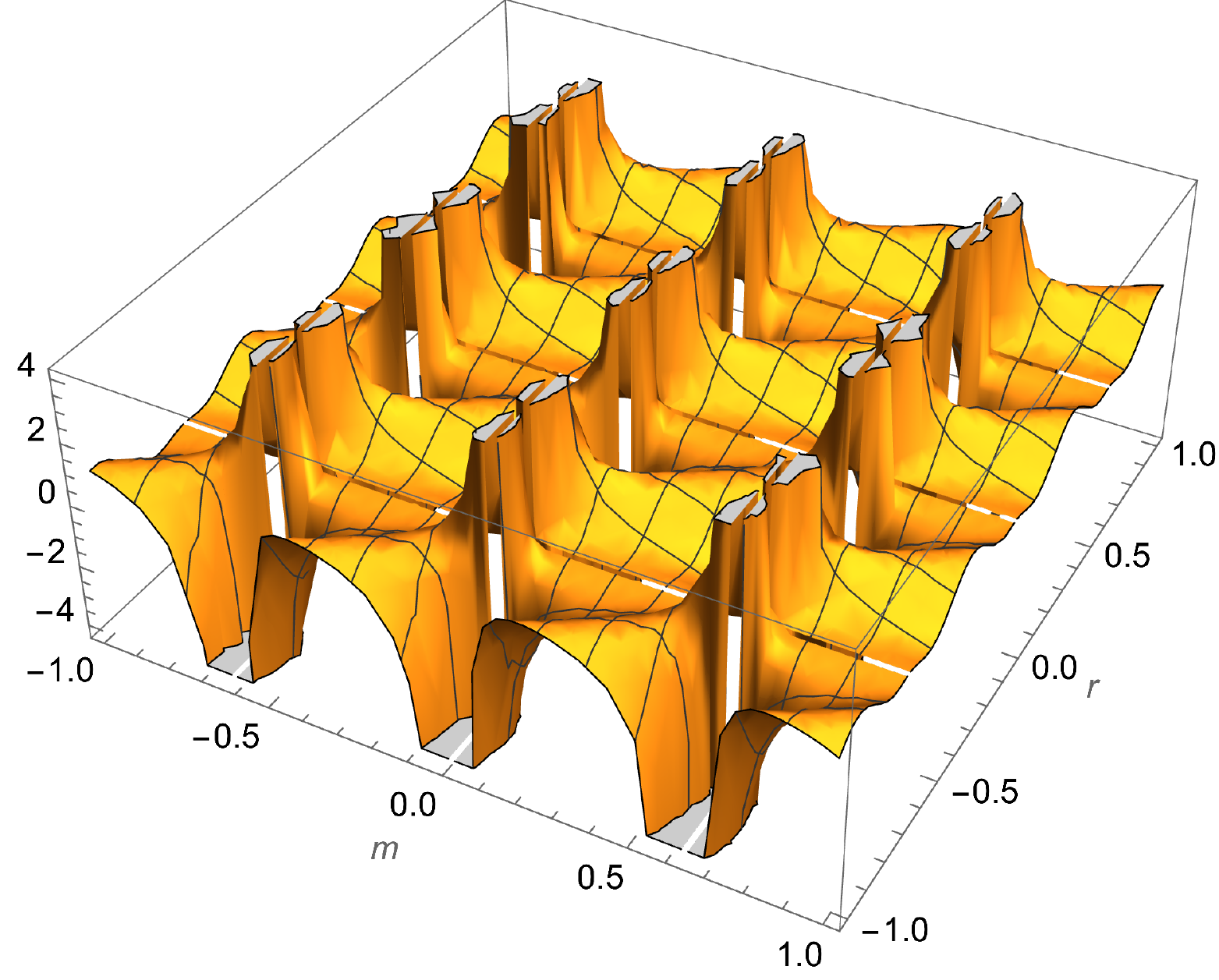}
\caption{$\text{Re}\left(\exp \left(2 e^{2 i m n} \Phi \left(e^{2 i m n},1,1+\frac{1}{2 n}\right)-2 e^{2 i n r} \Phi \left(e^{2 i n r},1,1+\frac{1}{2 n}\right)\right)\right)$\\ where $m,r\in\mathbb{R}$ and $n=5$}
   \label{fig:fig3}
\end{figure}
\vspace{-6pt}
\begin{figure}[H]
\includegraphics[scale=0.5]{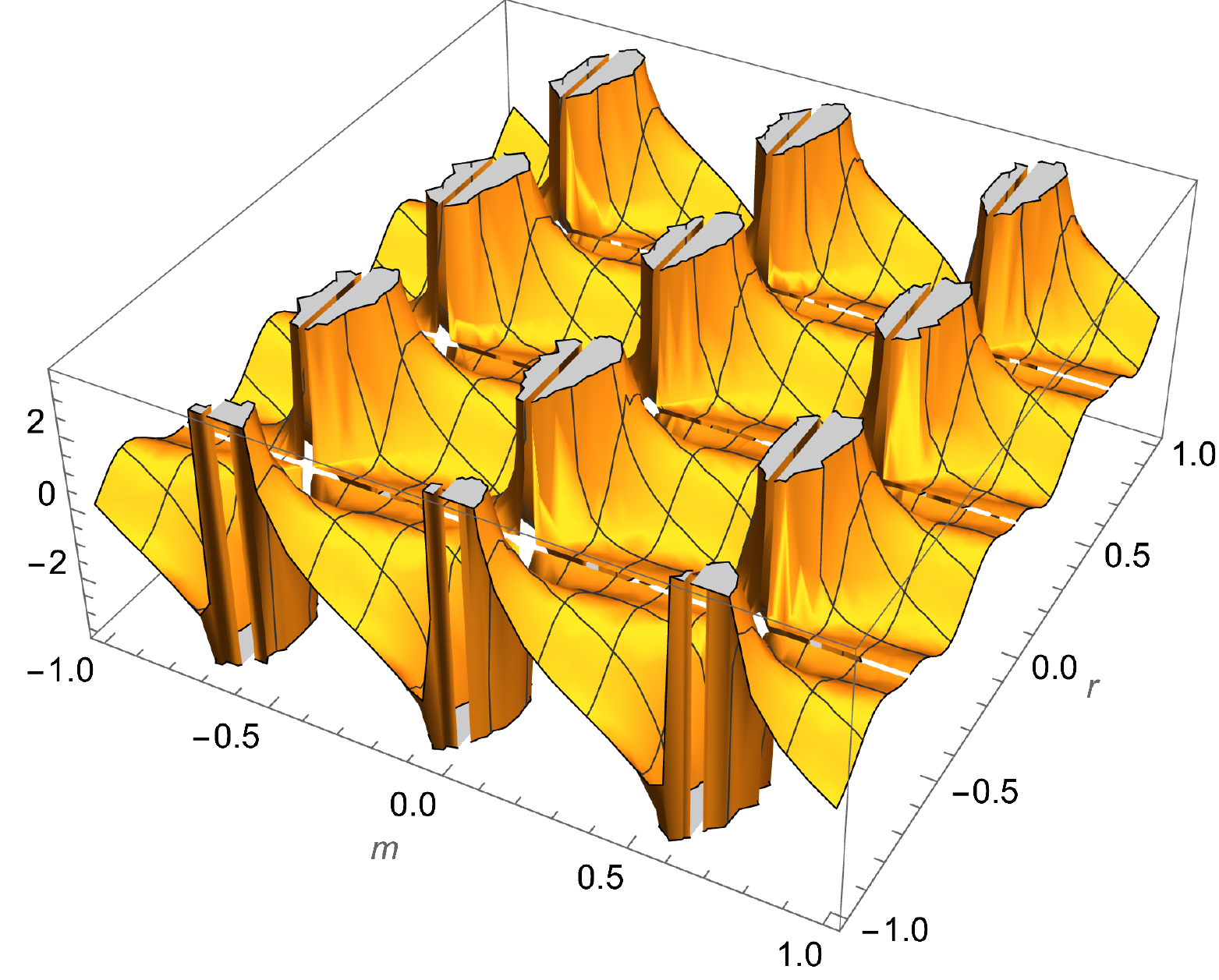}
\caption{$\text{Im}\left(\exp \left(2 e^{2 i m n} \Phi \left(e^{2 i m n},1,1+\frac{1}{2 n}\right)-2 e^{2 i n r} \Phi \left(e^{2 i n r},1,1+\frac{1}{2 n}\right)\right)\right)$.\\where $m,r\in\mathbb{R}$ and $n=5$}
   \label{fig:fig4}
\end{figure}
\vspace{-6pt}
\begin{figure}[H]
\includegraphics[scale=0.5]{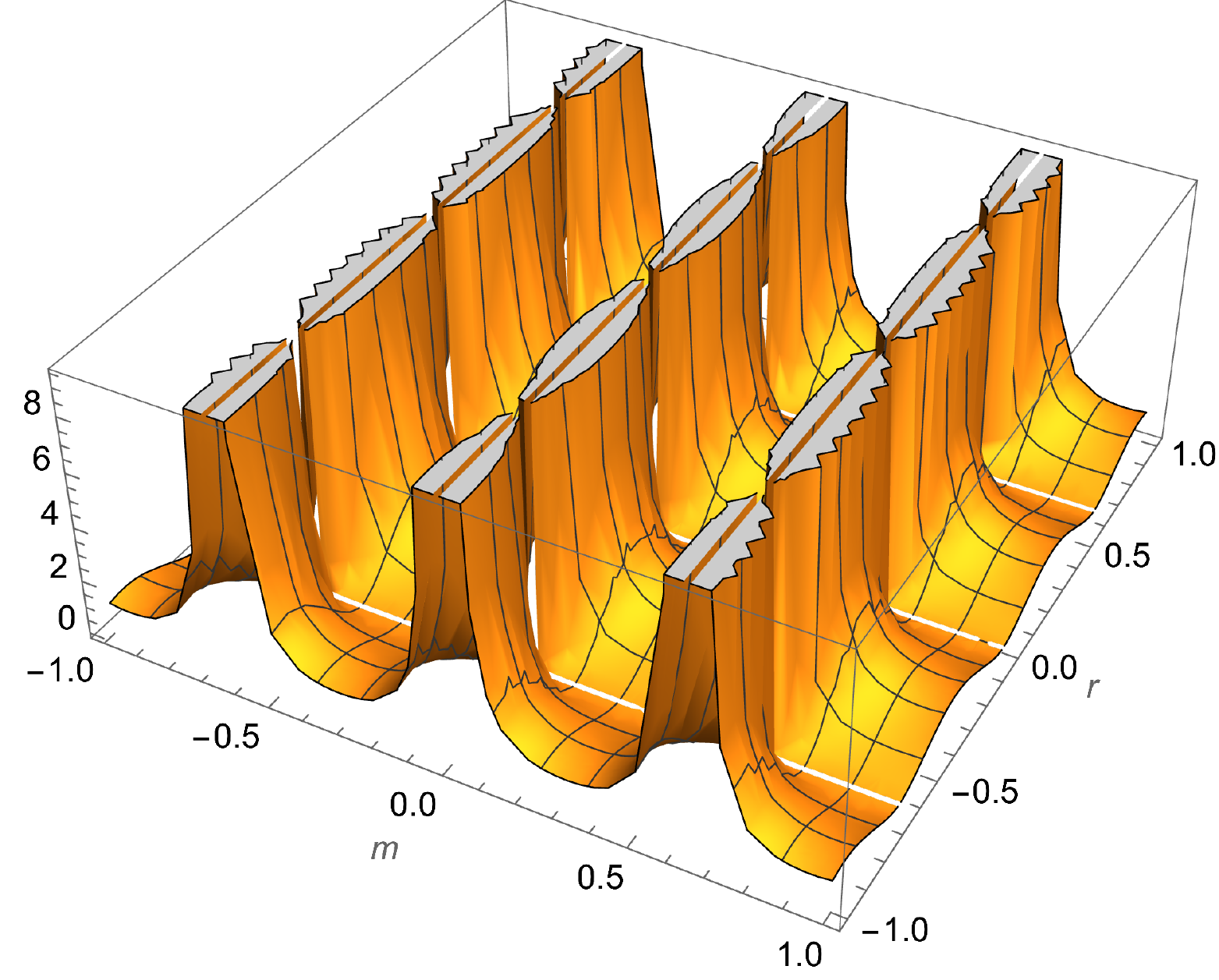}
\caption{$\text{Abs}\left(\exp \left(2 e^{2 i m n} \Phi \left(e^{2 i m n},1,1+\frac{1}{2 n}\right)-2 e^{2 i n r} \Phi \left(e^{2 i n r},1,1+\frac{1}{2 n}\right)\right)\right)$.\\where $m,r\in\mathbb{R}$ and $n=5$}
   \label{fig:fig4}
\end{figure}
\vspace{-6pt}
\begin{example}
Quotient Of Trigonometric Functions Raised To A Complex Power
\begin{multline}
\prod_{p=0}^{n-1}\frac{ \sin \left(\frac{l p \pi }{n}+\frac{x}{b}\right)}{\sin \left(\frac{l p \pi
   }{n}+x\right)}\left(\frac{\cos \left(\frac{l p \pi }{n}+\frac{x}{b}\right) \tan \left(\frac{l p \pi
   }{n}+\frac{x}{b}\right)}{\cos \left(\frac{l p \pi }{n}+\frac{x}{b^2}\right) \tan \left(\frac{l p \pi
   }{n}+\frac{x}{b^2}\right)}\right)^b\\
   =\frac{ \sin \left(\frac{n x}{b}\right)}{\sin (n
   x)}\left(\frac{\cos \left(\frac{n x}{b}\right) \tan \left(\frac{n x}{b}\right)}{\cos
   \left(\frac{n x}{b^2}\right) \tan \left(\frac{n x}{b^2}\right)}\right)^b
\end{multline}
\end{example}
\begin{proof}
Use equation (\ref{fslf}) setting $k=1,a=1,m=x,j=-1$ and simplify using the method in section (8.1) in \cite{reyn_ejpam}.
\end{proof}
\begin{example}
 Lerch Function Recurrence Identity With Consecutive Neighbours
\begin{multline}\label{eq:lerch1}
e^{i \pi  j l} \left(\Phi (z,s,a-j)+z^{2 j+1} \left(\Phi (z,s,a+j+1)+e^{i \pi  (j+1) l} \Phi \left(e^{i l \pi }
   z,s,a+j+1\right)\right)\right) \\
   +\Phi \left(e^{i l \pi } z,s,a-j\right)\\
   =2 e^{i \pi  j l} z^j \left(2^{1-s} z^2 \Phi
   \left(z^2,s,\frac{a}{2}+1\right)+e^{\frac{i \pi  s}{2}} (i a)^{-s}\right)
\end{multline}
\end{example}
\begin{proof}
Use equation (\ref{fslf}) and set $n=2,m=\log(z)/(2i),k=-s,a=e^{ai}$ and simplify.
\end{proof}
\begin{example}
Special Case Of A Recurrence Lerch Function Identity
\begin{equation}\label{eq:dlogp}
\Phi'(-i,0,a)+\Phi'(i,0,a)=2 \log \left(\frac{2 \Gamma
   \left(\frac{a}{4}\right)}{(a-2) \Gamma \left(\frac{a-2}{4}\right)}\right)
\end{equation}
\end{example}
\begin{proof}
Use equation (\ref{eq:lerch1}) and set $a=2a,l=1,z=i,j=-1$ and simplify using equation (25.14.2) in \cite{dlmf}. Next take the first partial derivative with respect to $s$ and set $s=0$ and simplify using equation (25.11.18 ) in \cite{dlmf}.
\end{proof}
\begin{example}
Special Case Of A Recurrence Lerch Function Identity
\begin{equation}\label{eq:dlogn}
\Phi'(-i,0,a)-\Phi'(i,0,a)=i \log \left(\frac{4 \Gamma
   \left(\frac{a+3}{4}\right)^2}{\Gamma \left(\frac{a+1}{4}\right)^2}\right)
\end{equation}
\end{example}
\begin{proof}
Use equation (\ref{eq:lerch1}) and set $a=2a,l=1,z=i,j=1$ and simplify using equation (25.14.2) in \cite{dlmf}. Next take the first partial derivative with respect to $s$ and set $s=0$ and simplify using equation (25.11.18 ) in \cite{dlmf}.
\end{proof}
\begin{example}
Special Case Of A Recurrence Lerch Function Identity
\begin{equation}
\Phi'(-i,0,a)^2-\Phi'(i,0,a)^2=2 i \log
   \left(\frac{\Gamma \left(\frac{a}{4}\right)}{2 \Gamma \left(\frac{a+2}{4}\right)}\right) \log \left(\frac{4 \Gamma
   \left(\frac{a+3}{4}\right)^2}{\Gamma \left(\frac{a+1}{4}\right)^2}\right)
\end{equation}
\end{example}
\begin{proof}
Use equations (\ref{eq:dlogp}) and (\ref{eq:dlogn}) and multiply to yield quoted result. 
\end{proof}
\begin{figure}[H]
\includegraphics[scale=0.5]{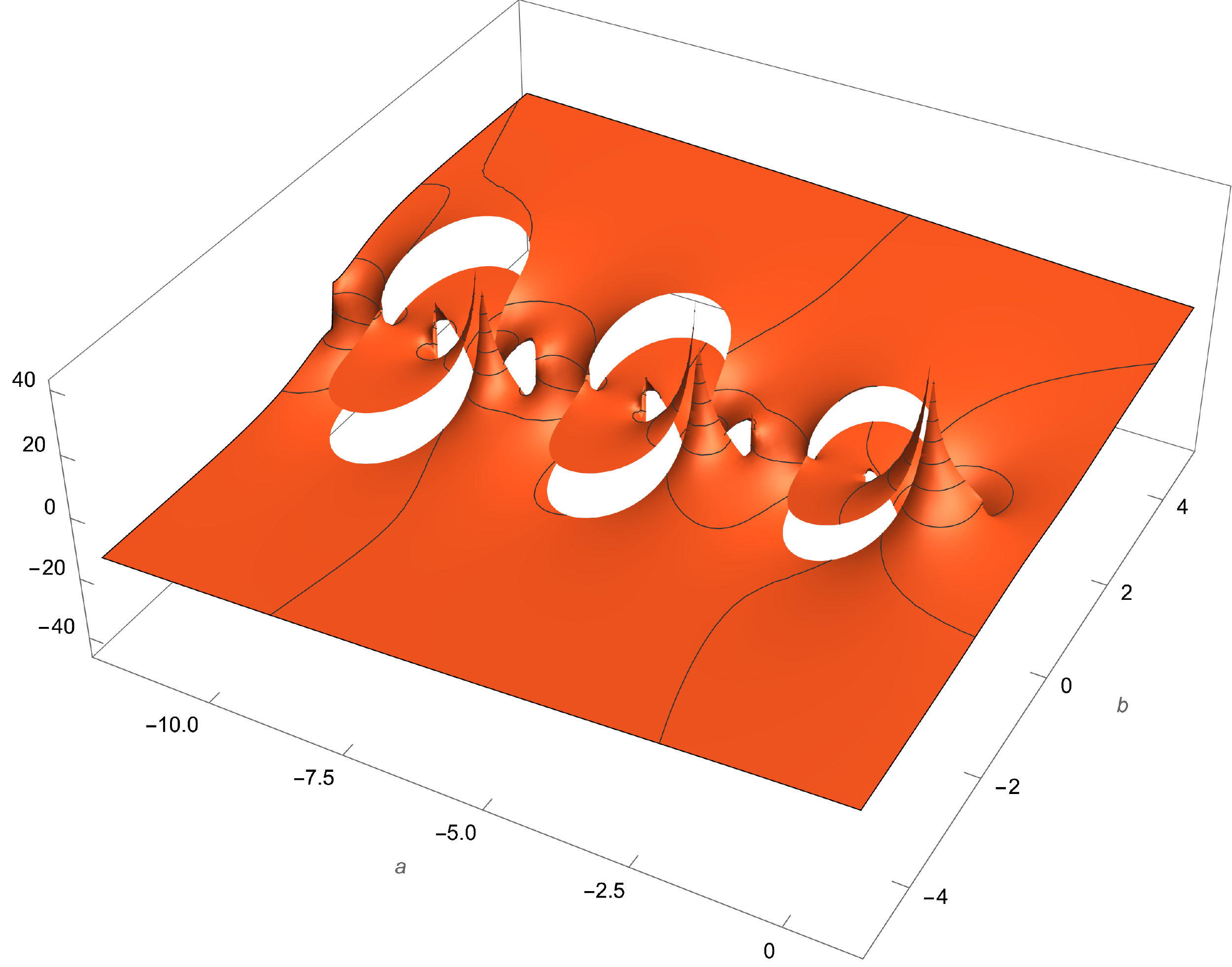}
\caption{$\text{Re}\left(2 i \log \left(\frac{\Gamma \left(\frac{1}{4} (a+i b)\right)}{2 \Gamma \left(\frac{1}{4} (a+i b+2)\right)}\right) \log
   \left(\frac{4 \Gamma \left(\frac{1}{4} (a+i b+3)\right)^2}{\Gamma \left(\frac{1}{4} (a+i b+1)\right)^2}\right)\right)$\\ where $a,b\in\mathbb{R}$}
   \label{fig:fig3}
\end{figure}
\vspace{-6pt}
\begin{figure}[H]
\includegraphics[scale=0.5]{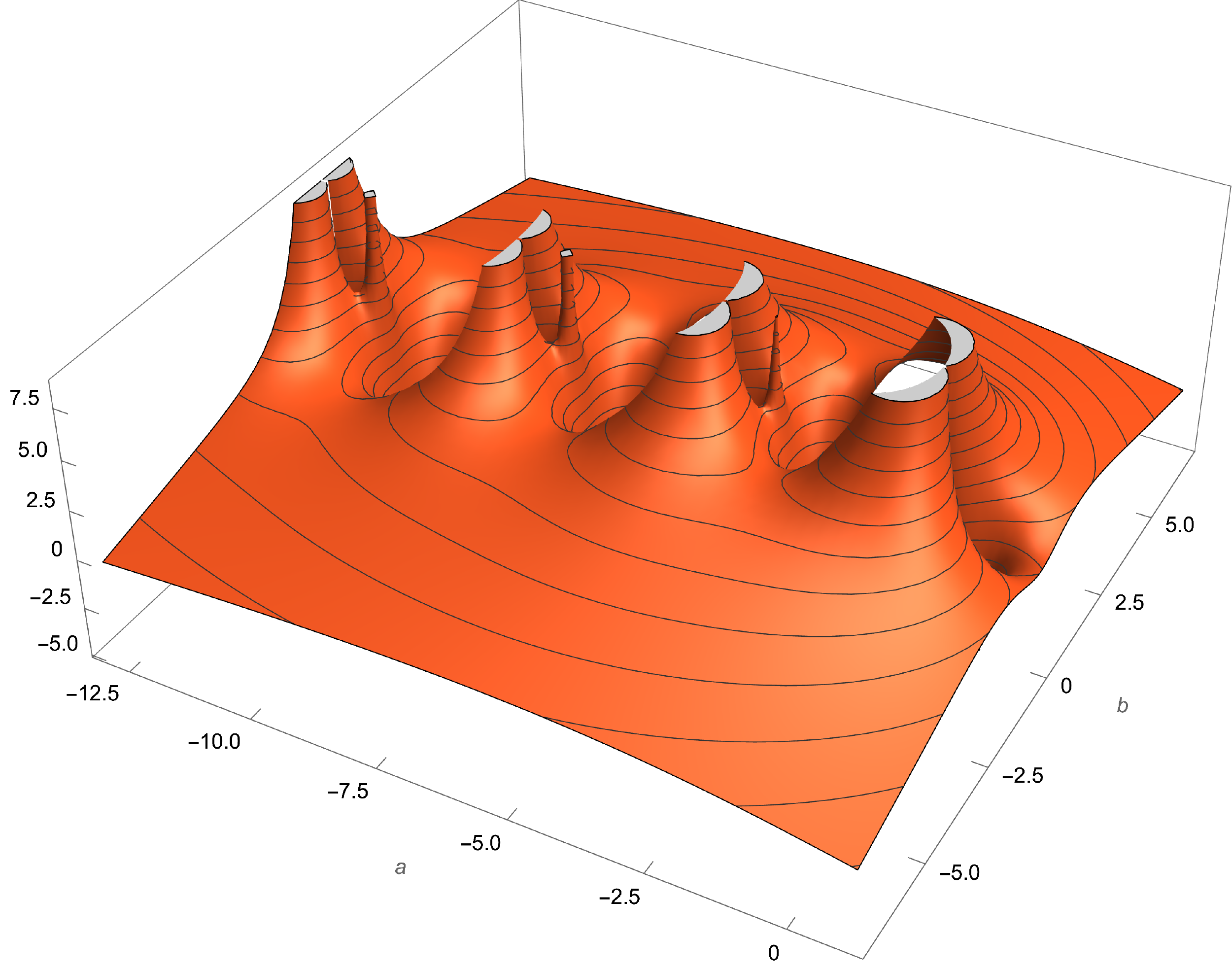}
\caption{$\text{Im}\left(2 i \log \left(\frac{\Gamma \left(\frac{1}{4} (a+i b)\right)}{2 \Gamma \left(\frac{1}{4} (a+i b+2)\right)}\right) \log
   \left(\frac{4 \Gamma \left(\frac{1}{4} (a+i b+3)\right)^2}{\Gamma \left(\frac{1}{4} (a+i b+1)\right)^2}\right)\right)$\\ where $a,b\in\mathbb{R}$}
   \label{fig:fig4}
\end{figure}
\vspace{-6pt}
\begin{figure}[H]
\includegraphics[scale=0.5]{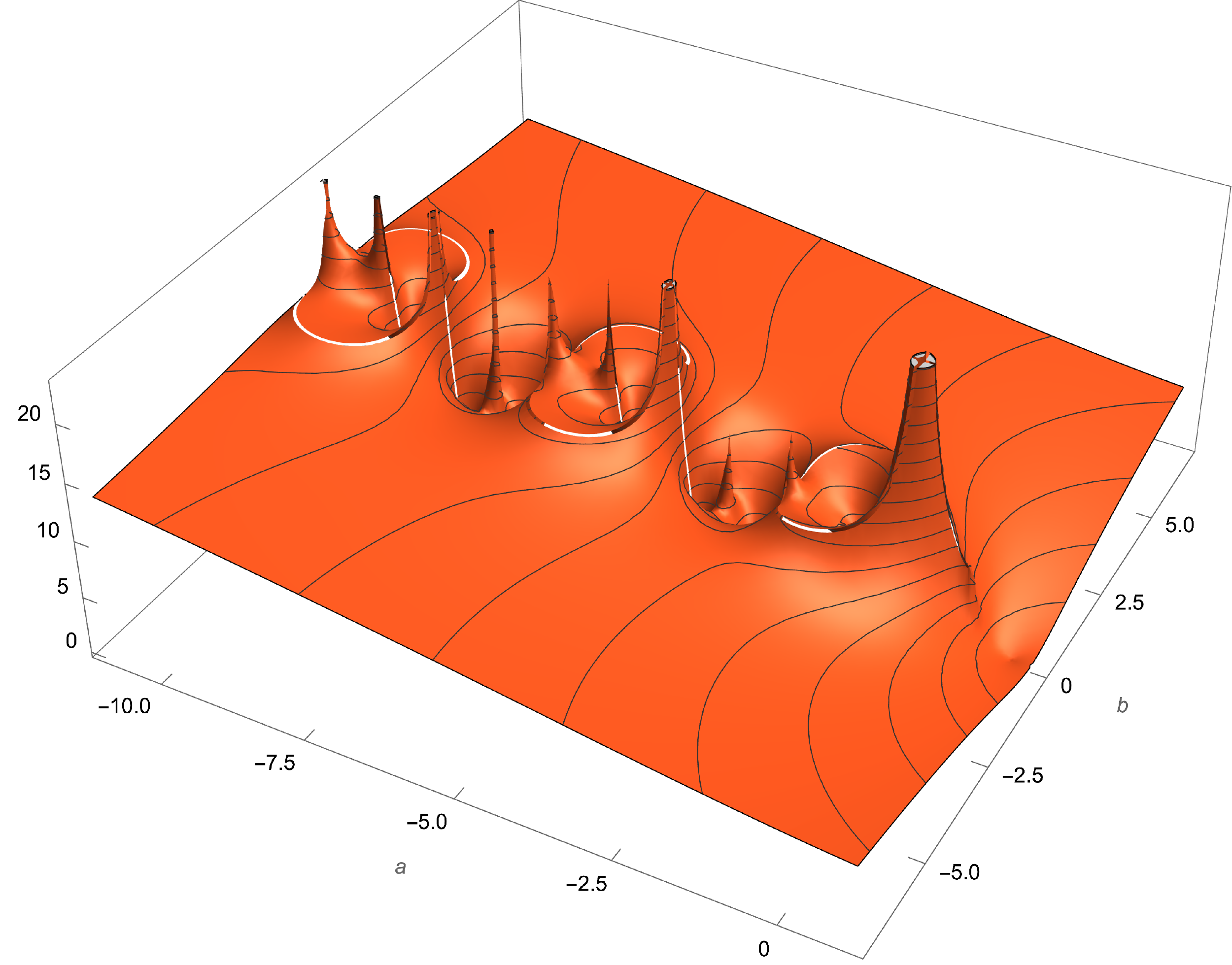}
\caption{$\text{Abs}\left(2 i \log \left(\frac{\Gamma \left(\frac{1}{4} (a+i b)\right)}{2 \Gamma \left(\frac{1}{4} (a+i b+2)\right)}\right) \log
   \left(\frac{4 \Gamma \left(\frac{1}{4} (a+i b+3)\right)^2}{\Gamma \left(\frac{1}{4} (a+i b+1)\right)^2}\right)\right)$\\ where $a,b\in\mathbb{R}$}
   \label{fig:fig4}
\end{figure}
\vspace{-6pt}
\section{Conclusion}
In this paper, we have presented a method for deriving a finite sum of the Hurwitz-Lerch Zeta function along with some interesting finite products of trigonometric functions using contour integration. We will be applying this method to other trigonometric formulae to derive other finite sums of the Hurwitz-Lerch function in future work. The results presented were numerically verified for both real and imaginary and complex values of the parameters in the integrals using Mathematica by Wolfram.
\end{document}